# HEEGAARD SPLITTINGS OF COMPACT 3-MANIFOLDS

MARTIN SCHARLEMANN

## Contents

Partially supported by NSF grant DMS 9504438

1. BACKGROUND

Here is a simple way to build a complicated 3-manifold. Begin with the 3-ball $B^3$ and in its boundary pick out two disjoint 2-disks $D_0$ and $D_1$. Using those disks, attach to $B^3$ a *handle*, that is a copy of $D^2 \times I$, by identifying $D^2 \times \{i\}$ with $D_i, i = 0, 1$. Depending on the orientation with which the ends of the handle are attached, the result is either $D^2 \times S^1$ or the non-orientable disk bundle over $S^1$, bounded by the Klein bottle. One can continue to add more handles to $B^3$ in a similar way. The result of adding $g$ of them is called a *genus $g$ handlebody*. Topologically, there are exactly two of them for any $g$, one of them orientable and the other not orientable, for once a non-orientable handle is attached, the end of any other handle can be slid over it, converting an orientable handle to a non-orientable, and vice versa. These manifolds are easily understood and not yet very complicated.

Now imagine taking two such handlebodies, $H_1$ and $H_2$, of the same genus and orientability. Then $\partial H_1$ and $\partial H_2$ are homeomorphic (the connected sum of $g$ tori or Klein bottles) and one can construct a complicated 3-manifold by attaching $H_1$ to $H_2$ by a possibly complicated homeomorphism of their boundaries. The resulting closed 3-manifold $M$ can be written $M = H_1 \cup_S H_2$, where $S$ is the surface $\partial H_i$ in $M$. This structure on $M$ is called a *Heegaard splitting* of $M$ and $S$ is a *splitting surface* (of a Heegaard splitting). Two Heegaard splittings of $M$ are *isotopic* if their splitting surfaces are isotopic in $M$. They are *homeomorphic* if there is a homeomorphism of $M$ carrying the splitting surface of one to the splitting surface of the other.

This method of constructing 3-manifolds is attributed to Heegaard [He] (see [Prz] for a translation into English of the relevant parts) though it was probably known to Poincare.

Natural questions arise: How *universal* is this construction? That is, how many closed 3-manifolds have such a structure? Is there a natural extension to 3-manifolds with boundary? This question is considered in section 2. How *unique* is such a structure? That is, given two such structures on the same 3-manifold, how are they related? This question is addressed in sections 6 and 7. How *useful* is the structure? That is, what information about the 3-manifold can be gleaned from the structure of a Heegaard splitting. Such questions are addressed in 5 and 8.

A useful earlier survey of the subject is [Zi], which focuses on Heegaard diagrams and on group presentations (briefly discussed in sections 2.3 and 5 below. I've relied heavily on its historical account. A central recent development has been an understanding of the importance of *strongly irreducible* Heegaard splittings (see 3.3), so their role has been chosen as a major theme of this survey.

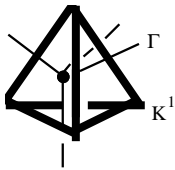

Figure 1.

## 2. Heegaard splittings and their guises

**2.1. Splittings from triangulations.** A foundational theorem of Moise [Mo] (see also [Bn]) says that all 3-manifolds can be triangulated. That is, given a compact connected 3-manifold $M$ there is a finite simplicial complex $K$ which is homeomorphic to $M$. For our purposes there are two important connected finite graphs in such a triangulation $K$: the 1-skeleton $K^1$ and the *dual* 1-skeleton $\Gamma$, defined as follows. (See fig. 1.) The vertices of $\Gamma$ are the barycenters of the 2- and 3-simplices of $K$ and edges connect the barycenter of a 3-simplex to the barycenter of each of its faces.

In case $M$ is closed, each 2-simplex is the face of precisely two 3-simplices, so each vertex in $\Gamma$ coming from a 2-simplex has valence 2. The edges of $\Gamma$ incident to such a vertex can therefore be amalgamated into a single edge, so that $\Gamma$ becomes a graph in which each vertex corresponds to a 3-simplex and each edge to the 2-simplex which it intersects. Now the regular neighborhood of a finite graph in a 3-manifold is easily seen to be a handlebody of genus $|edges| - |vertices| + 1$, for a regular neighborhood of a maximal tree is just a 3-ball, and a regular neighborhood of each remaining edge contributes a 1-handle. Furthermore, the region between regular neighborhoods of $K^1$ and $\Gamma$ is a product region, as can be seen easily in each 3-simplex. So, after thickening the regular neighborhoods of these graphs, $M$ can be viewed as the union of a regular neighborhood of $K^1$ and a regular neighborhood of $\Gamma$ along their common boundary. This is a Heegaard splitting of $M$.

**2.2. Splitting 3-manifolds with boundary.** The construction of Heegaard splittings for closed 3-manifolds in section 2.1 suggests several possible ways of extending the definition of Heegaard splitting to cover the case in which the 3-manifold has boundary. The most useful is the following: Write $\partial M$ as the disjoint union of two sets of components, $\partial_1 M$ and $\partial_2 M$. Choose a triangulation $K$ fine enough so that no simplex is incident to more than one boundary component. Let $K'$ be its barycentric subdivision. Delete the interior of all simplices of $K'$ incident to $\partial_2 M$. The resulting 3-manifold $M'$ is homeomorphic to $M$, since only a collar of $\partial_2 M$ has been removed; let $\partial'_2 M$ denote $\partial_2 M$ in this new triangulation. (Then $\partial'_2 M$ contains the subcomplex of $\Gamma$ determined by simplices incident to $\partial_2 M$.) Let $\Gamma_1 \subset M'$ be the union of $\partial_1 M$ and all vertices and edges of $K$ not incident to $\partial_2 M$. Let $\Gamma_2$ be the union of $\partial'_2 M$ and all vertices and edges of the dual 1-complex $\Gamma \cap M'$ not incident to $\Gamma_1$. Again it's easy to check that $M$ is the union of a regular neighborhood of the complexes $\Gamma_1$ and $\Gamma_2$ along their homeomorphic boundary, which is still a closed connected surface.

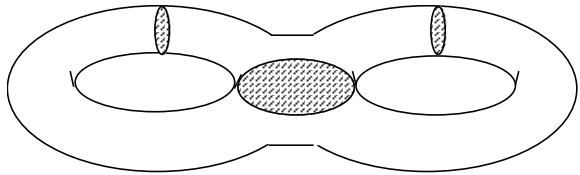

Figure 2.

This construction suggests the following way of defining a Heegaard splitting on a 3-manifold with boundary. A *compression body* $H$ is a connected 3-manifold obtained from a closed surface $\partial_- H$ by attaching 1-handles to $\partial_- H \times \{1\} \subset \partial_- H \times I$. (It is conventional to consider a handlebody to be a compression body in which $\partial_- H = \emptyset$.) Dually, a compression body is obtained from a connected surface $\partial_+ H$ by attaching 2-handles to $\partial_+ H \times \{1\} \subset \partial_+ H \times I$ and 3-handles to any 2-spheres thereby created. The cores of the 2-handles are called *meridian disks* and a collection of meridian disks is called *complete* if each of its complementary components is either a ball or $\partial_- H \times I$ (See fig. 2 for $H$ a handlebody). Suppose two compression bodies $H_1$ and $H_2$ have $\partial_+ H_1 \simeq \partial_+ H_2$. Then glue $H_1$ and $H_2$ together along $\partial_+ H_i = S$. The resulting compact 3-manifold $M$ can be written $M = H_1 \cup_S H_2$ and this structure also is called a *Heegaard splitting* of the 3-manifold with boundary $M$ (or, more specifically, of the triple $(M; \partial_- H_1, \partial_- H_2)$). It follows from the motivating discussion above that every compact 3-manifold has a Heegaard splitting.

2.3. **Splittings as handle decompositions - Heegaard diagrams.** Those familiar with handle decompositions of compact $n$-manifolds (see e. g. [RS, chapter 6] for the notation and viewpoint used here) will recognize similarities between the ways in which Heegaard splittings and handle decompositions are derived from a triangulation. The similarity goes deeper. Suppose $H_1 \cup_S H_2$ is a Heegaard splitting of a 3-manifold $(M; \partial_1 M, \partial_2 M)$. Then $H_1$ is obtained from $\partial_1 M \times I$ by attaching 1-handles and $H_2$ is obtained from $S = \partial_+ H_1 = \partial_+ H_2$ by attaching 2- and 3-handles. From this point of view a Heegaard splitting is just a standard handle decomposition of $M$ viewed as a cobordism between $\partial_1 M$ and $\partial_2 M$.

There is an advantage to this point of view. It is a standard trick in handle theory that the order of handles can frequently be rearranged. Always $r$-handles can be attached before $r+1$-handles, so that handles can be attached in ascending order. It is not generally true that an $(r+1)$-handle can be attached before an $r$-handle - it's necessary and sufficient that the attaching $r$-sphere of the $(r+1)$-handle be disjoint from the belt $(n-r-1)$-sphere of the $r$-handle. Translated into the language of Heegaard splittings this means that the natural order of handles can be rearranged if and only if there are essential disks in $H_1$ and $H_2$ whose boundaries are disjoint in $S$. This is a situation whose importance we will discuss later (see 3.3).

In this handle picture, all the topological information is contained in an understanding of the 1- and 2-handles, since the remaining 3-handles (if any) are uniquely determined by the spherical components of the boundary. Encouraged by this observation, we look for an efficient way of describing the way in which 2- handles are attached. We consider the case in which $M$ is closed; if $M$ has boundary the situation is analogous but a bit more complicated. When $M$ is closed, $H_1$ is a genus $g$ handlebody. The attaching curves $\partial \Delta_2$ for the cores $\Delta_2$ of the 2-handles constitute a family of simple closed curves in $\partial H_1$. We may as well isotope $\partial \Delta_2$ to intersect a chosen minimal complete collection $\Delta_1$ of meridian disks for $H_1$ transversally and minimally. When $H_1$ is cut open along $\Delta_1$ it becomes a 3-ball $B^3$, on whose boundary appear two copies of each disk of $\Delta_1$. Let $V \subset \partial B^3$ be this collection of disks. The attaching curves $\partial \Delta_2$ are (typically) also cut up - into a collection $\mathcal{A}$ of arcs and simple closed curves in $\partial B^3 - V$. The ends of each arc in $\mathcal{A}$ lie in $\partial V$.

If the splitting is irreducible (see (see 3.2), note that $\mathcal{A}$ consists entirely of arcs, since any simple closed curve in $\mathcal{A}$ bounds a disk in $B^3$ and the union of this disk and a 2-handle core with the same boundary would give a reducing sphere. When no component of $\mathcal{A}$ is a simple closed curve, we can think of the union of $V$ and $\mathcal{A}$ as defining a graph $\Gamma$ in $\partial B^3$, with fat vertices $V$ and edges $\mathcal{A}$. There is additional structure, of course, which identifies each pair of vertices of $V$ that began as the same disk in $\Delta_1$ and also identifies ends of edges that were cut at $\partial \Delta_1$.

The graph $\Gamma$ has some pleasant properties. For example, there are no trivial loops in $\Gamma$, for such a loop could have been removed by an isotopy of $\partial \Delta_2$ that lowers $\partial \Delta_1 \cap \partial \Delta_2$. But there are a lot of choices made in the construction or $\Gamma$ (e. g. $\Delta_1$ and $\Delta_2$) so it is not particularly well-defined. The use of these diagrams to study the underlying 3-manifold can be quite complicated and is often disappointing.

Sometimes the ability to rechoose $\Delta_1$ and $\Delta_2$ can be useful. For example, although we have observed that $\Gamma$ contains no trivial loops, it is also true, when the splitting is irreducible, that if any loop at all appears, the diagram can be simplified. A loop in a Heegaard diagram (i. e. an edge in $\mathcal{A}$ both of whose ends lie on the same vertex in $V$) is sometimes called a *wave*. A wave, and the vertex $v$ in $V$ on which it is based, divides $\partial B^3$ into two disks. One of them (call it $E$) does not contain the other vertex in $V$ that is identified with $v$ in $\Delta_1$. All ends of 1-handles in $H_1$ represented by vertices that lie in $E$ can be dragged over the 1-handle in $H_1$ whose cocore is $v$ (thereby redefining $\Delta_1$). At this point the wave becomes an inessential loop, which can be removed by an isotopy. The net effect is to reduce $\Delta_1 \cap \Delta_2$ by redefining $\Delta_1$.

We refer the reader to the excellent [Zi] for a more thorough discussion of Heegaard diagrams.

2.4. **Splittings as Morse functions and as sweep-outs.** Smooth manifolds admit handle structures just as $PL$ manifolds do. One way of showing this classical fact is via Morse theory [Mi]. A generic smooth height function $h$ from the smooth manifold $M$ to $R$ will have only non-degenerate critical points. At each critical height $t_0$, as the part $h^{-1}(-\infty, t_0 - \epsilon]$ of $M$ lying below $t_0$ changes to $h^{-1}(-\infty, t_0 + \epsilon]$, the topological effect is to add a handle. The handles can then be rearranged to appear in ascending order, just as in the PL theory.

The argument is reversible. Given a handle structure on a smooth manifold $M$ one can easily construct a Morse function which induces that handle structure. So associated to a Heegaard splitting there is also a Morse function. That is, given a Heegaard structure $H_1 \cup_S H_2$ on $(M; \partial_1 M, \partial_2 M)$ there is a Morse function $h : M \to [0, 1]$ with singular values $0 < a_1 < a_2 < \ldots < a_k < b_1 < \ldots b_j < 1$ so that

1. $a_i$ is an index one critical level
2. $b_j$ is an index two critical level
3. $h^{-1}(0) = \partial_1 M$ or is an index zero critical point if $\partial_1 M = \emptyset$
4. $h^{-1}(1) = \partial_2 M$ or is an index three critical point if $\partial_2 M = \emptyset$
5. if $a_k < t < b_1$ then $h^{-1}(t) \cong S$.

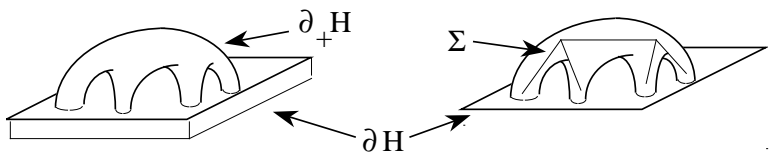

Figure 3.

A reason for taking this viewpoint is that it is sometimes advantagous to put knots and graphs in $M$ into the simplest possible position with respect to the Heegaard splitting. One way to accomplish this is to incorporate Gabai's powerful notion of "thin position" into the theory and make the knot or graph thin with respect to this Morse function (cf. [BO1], [ST1], 3.7, 4.1).

A similar way to use a Heegaard splitting to parameterize the 3-manifold $M$ is to focus attention on heights between the top index-one critical level $a_k$ and the bottom index-two critical level $b_1$, as we now explain.

Define the *spine* $\Sigma$ of a handlebody $H$ to be a finite graph in $H$ for which $H$ is a regular neighborhood. From the construction, every handlebody has a spine (for a closed triangulated manifold, the 1-skeleton and the dual 1-skeletons are spines of the relevant handlebodies in the construction described in section 2.1 above). The spine of $H$ is not uniquely defined, but any two spines differ by a sequence of "edge-slides" (see [ST1, 1.2]). For $H$ a compression body, a spine $\Sigma$ is a graph in $H$ so that $\Sigma \cap \partial H = \Sigma \cap \partial_- H$ consists only of valence one vertices and $H$ deformation retracts to $\Sigma \cup \partial_- H$. (See fig. 3.) Again the construction of $H$ guarantees the existence of a spine, and two spines of the same compression body differ by a series of edge slides, where ends of edges may be slid along paths in $\partial_- H$.

Notice that the complement of a spine in $H$ is homeomorphic to $\partial_+ H \times I$. Suppose then we are given a Heegaard splitting $H_1 \cup_S H_2$ of $M$ and spines $\Sigma_i$ of each $H_i$. Then $M - (\Sigma_1 \cup \Sigma_2)$ is just a product $S \times I$. This parameterization of $M - (\Sigma_1 \cup \Sigma_2)$ is sometimes called a "sweep-out" by $S$ since $S$ sweeps between one spine and the other. This viewpoint allows great flexibility in the positioning of the splitting surface. (See section 7.6).

A related idea is to consider a single spine, say $\Sigma_1 \subset M$, as a graph in $M$ for which there are $\partial$-singular compressing disks (the meridian disks of $\Sigma_2$). In some situations the $\partial$-singular disks can be used to slide $\Sigma_1$ into useful positions. (See [Ot] and [ST2].)

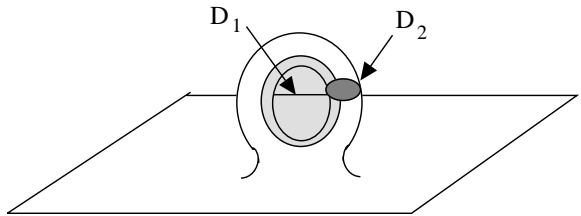

Figure 4.

## 3. Structures on Heegaard splittings

3.1. **Stabilization.** As we have seen, Heegaard splittings have connections to triangulations, handle decompositions, and Morse functions on 3-manifolds. Just as triangulations can be subdivided, or a Morse function locally perturbed to introduce cancelling critical points, or, in a handlebody description, a cancelling pair of handles can be inserted, so there is a natural and trivial way of making a Heegaard splitting more complicated. Suppose $H_1 \cup_S H_2$ is a Heegaard splitting of a 3-manifold $M$ and $\alpha$ is a properly imbedded arc in $H_2$ parallel to an arc in $S$. Here "parallel" means that there is an embedded disk $D$ in $H_2$ whose boundary is the union of $\alpha$ and an arc in $\partial_+ H_2$. Now add a neighborhood of $\alpha$ to $H_1$ and delete it from $H_2$. This adds a 1-handle to $H_1$ (whose core is $\alpha$) and, topologically, also adds a 1-handle to $H_2$ (whose cocore is $D$). So once again the result is a Heegaard splitting $H_1' \cup_{S'} H_2'$, where the genus of each $H_i'$ is one greater than $H_i$. This process is called a *stabilization* of $S$.

Stabilization is uniquely defined. That is, any two splittings obtained by stabilizing the same splitting surface are isotopic. On the other hand, there is no reason to believe that if the stabilizations of two different splitting surfaces are isotopic, then the original surfaces were isotopic. Indeed, it will turn out (see Section 7) that any two Heegaard splittings of the same manifold will become isotopic after a sufficient number of stabilizations. So in stabilizing a splitting surface we may lose and can't gain information about its structure. Interest therefore focuses on splittings which are not stabilizations of other splittings, that is splittings which cannot be destabilized. How is this detectable? (See fig. 4.)

**Lemma 3.1.** *A splitting $M = H_1 \cup_S H_2$ can be destabilized if and only if there are properly imbedded disks $D_i \subset H_i$ so that $|\partial D_1 \cap \partial D_2| = 1$.*

**Proof**: Suppose a splitting is stabilized as above. Then let $D_1$ be the cocore disk of the 1-handle attached along $\alpha$ and let $D_2$ be the disk $D$.

Conversely, suppose disks $D_i \subset H_i$ are as in the lemma. Because the boundaries intersect in a single point, each disk is non-separating and hence essential. Let $T_i$ be the surface obtained from $S$ by compressing along $D_i$, converting $H_i$ into a simpler compression body $J_i$. The union of a bicollar of $D_1$ in $H_1$ and $D_2$ in $H_2$ along the square in which they intersect is a 3-ball intersecting each of $T_i$ in a hemisphere, and so defines an isotopy between the $T_i$. In particular $T_1$ divides $M$ into $J_1$ and an isotope of $J_2$ and so is a Heegaard splitting surface. It's easy to see that stabilizing $T_1$ gives $S$: the 1-handle dual to $D_1$ corresponds to the arc $\alpha$ and the disk $D_2$ corresponds to the disk $D$. □

### 3.2. Reducible splittings.

Suppose $M$ and $M'$ are two 3-manifolds with Heegaard splittings $H_1 \cup_S H_2$ and $H_1' \cup_{S'} H_2'$. From these we can naturally construct a Heegaard splitting of the connected sum $M" = M \# M'$ as follows: Remove from $M$ and $M'$ 3-balls $B$ and $B'$ which intersect $S$ and $S'$ respectively in equatorial disks. Glue together the boundaries of the 3-balls so that each hemisphere $H_i \cap \partial B$ is attached to $H_i' \cap \partial B'$. The resulting surface $S \# S'$ splits $M \# M'$ into compression bodies $H_i"$. To see that the complementary pieces are compression bodies, note that topologically $H_i"$ is obtained from the disjoint union of $H_i$ and $H_i'$ by attaching a 1-handle whose two ends lie in $\partial_+ H_i$ and $\partial_+ H_i'$ respectively.

Conversely, given a Heegaard splitting $H_1 \cup_S H_2$ of a 3-manifold $M"$ and a 2-sphere $P$ which intersects $M$ in a single circle, we can get a Heegaard splitting of the reduced 3-manifold, obtained by doing surgery on $P$ (the reduced manifold is the disjoint union of $M$ and $M'$ when $P$ is separating, as in the above example). If $P$ bounds a ball in $M"$ and $S$ intersects the ball in a single equatorial disk, then the manifolds $M"$ and $M$, say, are the same and get the same splitting. Otherwise, the splitting of the reduced manifold is simpler, since the genus of the splitting surface is reduced, and the Heegaard splitting of $M"$ can be easily reconstructed from the splitting of the reduced manifold. These considerations lead to the following:

**Definition 3.2.** *A Heegaard splitting $H_1 \cup_S H_2$ is* reducible *if there is a 2-sphere which intersects $S$ in a single essential circle.*

An alternate way of saying this is that there are essential properly imbedded disks $D_i \subset H_i$ so that $\partial D_1 = \partial D_2$ in $S$. (See fig. 5.)

There is a connection with stabilization, given by the following lemma:

**Proposition 3.3.** *Suppose $H_1 \cup_S H_2$ is a splitting that can be destabilized. Then either it is reducible or it is the standard genus one splitting of $S^3$.*

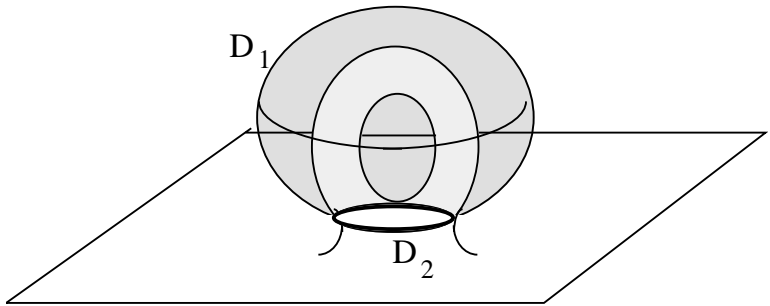

Figure 5.

**Proof**: Let $D_i \subset H_i$ be disks so that $\partial D_1 \cap \partial D_2$ is a single point. As in the proof of 3.1, let $B$ be the union of a bicollar of $D_1$ in $H_1$ and a bicollar of $D_2$ in $H_2$ along the square in which they intersect. $B$ is a 3-ball whose boundary sphere $P$ can be moved slightly (e. g. increase the radius of both disks slightly) so that $P$ intersects each $H_i$ in a single hemisphere and so that the curve $c = P \cap S$ (the boundary of a regular neighborhood of the figure-eight $\partial D_1 \cup \partial D_2$) cuts off from $S$ a punctured torus. Unless this curve is inessential in $S$ the boundary of the 3-ball is a reducing sphere. If the curve is inessential, then $S$ is a torus dividing $M$ into two solid tori, whose meridians intersect in a single point. This is the genus one Heegaard splitting of $S^3$. □

One of the first major theorems on Heegaard splittings, due to Haken, is that any Heegaard splitting of a reducible manifold is a reducible splitting. The theorem is important not just for what it says, but for the type of argument which is used.

**Theorem 3.4** ([Ha1]). *Suppose $M$ is a reducible manifold with a Heegaard splitting $H_1 \cup_S H_2$. Then there is a reducing sphere $P$ for $M$ so that $P \cap S$ is a single circle.*

**Proof**: There are two ways to do this. Similar to Haken's original proof is that given in [Ja, II.7]. One can assume that $P$ intersects (either) one of the compression bodies only in disks. (One way to do this is to put a spine of $H_2$, say, transverse to $P$ and take $H_2$ to be very thin.) The idea will be to minimize the number of circles of intersection, under the assumption that $P$ intersects one of the compression bodies only in disks. If $P$ intersects $H_2$ only in disks, consider the planar surface $P_1 = P \cap H_1$. Compress and $\partial$-compress $P_1$ as much as possible. Compressions of $P_1$ will convert $P$ into two spheres, at least one of which is a reducing sphere - restrict attention to that one. At the end of this process $P_1$ will be converted to a surface $P'$ which is disjoint from a complete collection of meridian disks for $H_1$ (otherwise curves of intersection can be used to compress or $\partial$-compress) and, for any essential curve $\alpha$ in $\partial_- H_1$, disjoint from a spanning annulus $\alpha \times I$ (same argument). It follows that $P' \cap H_1$ is a collection of disks. What is not obvious, but can be explicitly calculated, is that the number of disks in $P' \cap H_1$ is lower than the original $P \cap H_2$. The process is continued, switching the roles of $H_1$ and $H_2$ until there is only one intersection curve.

Another approach is given in [ST2]. Put $\Sigma_2$, the spine of $H_2$, transverse to $P$. Let $\Delta$ be a complete collection of compressing disks for $H_1$ viewed as a $\partial$-singular collection of disks in the complement of $\Sigma_2$. Put $\Delta$ transverse to $P$. Circles of intersection can be removed, just as in the previous argument, so that $(\Sigma_2 \cup \Delta) \cap P$ becomes a graph $\Gamma \subset P$ with vertices $\Sigma_2 \cap P$ and edges $\Delta \cap P$. Trivial loops of $\Gamma$ can be eliminated at the cost of merely changing $\Delta$, and a vertex incident to some edges but no loops can be used to slide edges of $\Sigma_2$ in a way that lowers $\Sigma_2 \cap P$. (This is the hard part to see.) The upshot is that, eventually, there is guaranteed to be an isolated vertex. This picks out a meridian $\mu$ of $H_1$ which is disjoint from a complete collection of meridian disks for $H_2$. If $H_2$ is a handlebody this implies that $\partial \mu$ also bounds a meridian in $H_2$ and so $H_1 \cup_S H_2$ is reducible. If $H_2$ is merely a compression body, we can only conclude that there is a $\partial$-reducing disk for $M$ which intersects $S$ in a single curve. But we can surger $M$ along this disk to get a new reducible 3-manifold and continue the process until an appropriate sphere is found. $\square$

The last step of the second proof suggests a new notion:

**Definition 3.5.** *A Heegaard splitting $M = H_1 \cup_S H_2$ is $\partial$-reducible if there is a $\partial$-reducing disk for $M$ which intersects $S$ in a single curve.*

It also suggests the following analogue to Theorem 3.4.

**Proposition 3.6.** *Any Heegaard splitting of a $\partial$-reducible 3-manifold is $\partial$-reducible.*

**Proof**: Both proofs above easily generalize. $\square$

A more difficult theorem, discussed in more detail in section 6.1 but relevant here, characterizes Heegaard splittings of the 3-sphere.

**Theorem 3.7** ([Wa]). *Every positive genus Heegaard splitting of $S^3$ is stabilized.*

This implies, more fully, that any positive genus Heegaard splitting of $S^3$ is obtained by stabilizing the unique genus zero splitting into 3-balls. So a Heegaard splitting of $S^3$ is completely determined by its genus.

Armed with Theorem 3.7 we can prove a sort of converse to 3.3

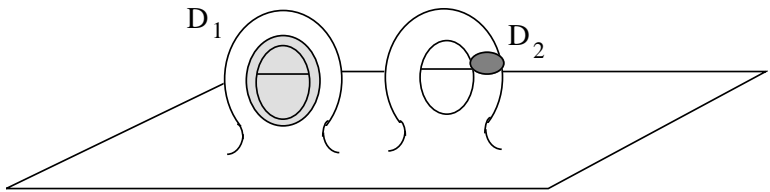

Figure 6.

**Theorem 3.8.** *Suppose $M$ is an irreducible 3-manifold and $H_1 \cup_S H_2$ is a reducible Heegaard splitting of $M$. Then $H_1 \cup_S H_2$ is stabilized.*

**Proof**: Let $P$ be a sphere which intersects $S$ in a single essential circle. Since $M$ is irreducible, $P$ bounds a 3-ball in $M$, so the manifold obtained by reducing $M$ along $P$ is the disjoint union of $S^3$ and a homeomorph of $M$. The induced Heegaard splitting of the former is, by 3.7, stabilized. Its stabilizing disks, when viewed back in $H_1 \cup_S H_2$ show that $S$ was also stabilized. □

3.3. **Weakly reducible splittings.** In 1987 ([CG]) Casson and Gordon discovered a new structure on Heegaard splittings which is perhaps less natural than those described above but which has turned out to be quite useful.

**Definition 3.9.** *A Heegaard splitting $H_1 \cup_S H_2$ is* weakly reducible *if there are essential disks $D_i \subset H_i$ so that $\partial D_1$ and $\partial D_2$ are disjoint in $S$.*

**Remarks:**
1. This notion coincides precisely to the assertion that, in viewing the Heegaard structure as a handle decomposition, at least one 2-handle ($D_2$) can be attached before all 1-handles are attached (in particular the 1-handle dual to $D_1$).
2. Any reducible Heegaard splitting is weakly reducible, simply by cutting the sphere $P$ that intersects $S$ into two disks along $P \cap S$, then pushing the two boundaries apart.
3. A splitting that is not weakly reducible is called *strongly irreducible*.

Here are two sample applications of this structure:

**Lemma 3.10** ([ST2])**.** *Suppose $H_1 \cup_S H_2$ is a strongly irreducible splitting of a 3-manifold $M$ and $F$ is a disk in $M$ transverse to $S$ with $\partial F \subset S$. Then $\partial F$ also bounds a disk in some $H_i$.*

**Proof**: The proof is by induction on $|S \cap int(F)|$. If the interior of $F$ is disjoint from $S$ there is nothing to prove. If $S - F$ has any disk components $D$ then, by replacing the subdisk of $F$ bounded by $\partial D$ by a parallel copy of $D$ we can decrease $|S \cap int(F)|$. So assume that each curve in $S \cap F$ is essential in $S$.

A disk component of $F - S$ compresses $S$ in one of the two compression bodies, say $H_1$. Then by strong irreducibility of $S$, all disk components of $F - S$ lie in $H_1$. If any pair of curves of $F \cap S$ are nested and inessential in $F$ then the outer curve of the innermost such pair cuts off a component $P$ of $F - S$ so that all but one of the curves in $\partial P$ are adjacent to disks in $H_1$ (hence $P \subset H_2$) and precisely one, denoted $\alpha$, is not. Compress $S$ into $H_1$ along 2-handles whose cores are the disks with boundaries on $\partial P$. Let $M_-$ be the 3-manifold obtained from $H_2$ by attaching these 2-handles to $H_2$. Then $\alpha \subset \partial M_-$ is inessential in $M_-$ so, by strong irreducibility and 3.6, $\alpha$ is inessential in $\partial M_-$. Push the disk $\alpha$ bounds in $\partial M_-$ slightly into $H_1$ and observe that this is then a disk $D$ in $H_1$ whose boundary is parallel to $\alpha$ in teh component of $F$ adjacent to $P$ across $\alpha$. Replacing the subdisk of $F$ bounded by $\alpha$ (or all of $F$ if $\alpha = \partial F$) with $D$ lowers $|S \cap int(F)|$. □

**Theorem 3.11** ([CG]). *If $M = H_1 \cup_S H_2$ is a weakly reducible splitting then either $H_1 \cup_S H_2$ is reducible or $M$ contains an incompressible surface.*

**Proof**: $S$ can be compressed simultaneously in both directions, that is, both into $H_1$ and simultaneously into $H_2$. Let $\Delta_1 \subset H_1$ and $\Delta_2 \subset H_2$ be collections of disjoint meridians in the respective compression bodies so that $\partial \Delta_1$ and $\partial \Delta_2$ are disjoint in $S$ and the families $\Delta_i$ are maximal with respect to this property. That is, if $S_i$ represents the surface in $H_i$ obtained by compressing $S$ along $\Delta_i$, then any further compressing disks of $S_i$ into $H_i$ will necessarily have boundaries intersecting the boundaries of the other disk family (or any obtained from it by 2-handle slides - a requirement that makes the definition of "maximal" here mildly subtle).

Let $\bar{S}$ be the surface obained by compressing $S_1$ along $\Delta_2$ (or, symmetrically, $S_2$ along $\Delta_1$). $\bar{S}$ separates $M$ into the remnant $W_1$ of $H_1$ and the remnant $W_2$ of $H_2$. Dually, $H_1$ can be recovered from $W_1$ by removing some "tunnels" (neighborhoods of arcs) from $W_1$ and attaching some 1-handles in $W_2$. A helpful and vivid picture is to imagine $H_1$ red and $H_2$ blue. The compressions of $S$ to $\bar{S}$ along the $\Delta_i$ cover $\bar{S}$ with both red and blue spots, two red spots for each disk in $\Delta_1$ and two blue spots for each disk in $\Delta_2$. $S$ is recovered from $\bar{S}$ by attaching red tubes in $W_2$ with ends on red spots and blue tubes in $W_1$ with ends on blue spots.

The surface $\bar{S}$ is incompressible in $M$. To see this, suppose that $\bar{S}$ compresses into $W_1$, say. After pushing $\bar{S}$ slightly into $W_2$, we can view $S_1$ as a Heegaard splitting surface of $W_1$, that is $W_1 = H_1 \cup_{S_1} (W_1 \cap H_2)$. The compression of $\bar{S}$ is a $\partial$-reduction of $W_1$. By Theorem 3.6 there is a $\partial$-reducing disk $D$ that intersects $S_1$ in a single circle. We can take $\partial D$ to be disjoint from the "red spots" (i. e. disjoint from $\Delta_1$) and, after some 2-handle slides among the $\Delta_2$, we can make $\Delta_2$ disjoint from the annulus $D - H_1$. But then $D \cap H_1$ makes $S_1$ compressible in $H_1$ via a disk disjoint from $\Delta_2$, contradicting the maximality of $\Delta_1$.

Unless $\bar{S}$ is a collection of spheres, we are through. Suppose $\bar{S}$ is a collection of spheres. Note that at least one, $\bar{S}_0$, has both a red spot and a blue spot. For otherwise, when $S$ is recovered from $\bar{S}$ by attaching red and blue tubes, $S$ would consist of two components: one containing all red tubes and one containing all blue. Choose in $\bar{S}_0$ a simple closed curve that separates in the sphere $\bar{S}_0$ the red spots from the blue spots. Push the interior of the disk in $\bar{S}_0$ that contains the red spots (resp. blue spots) completely into $H_1$ (resp. $H_2$). Then $\bar{S}_0$ is the union of a red disk and a blue disk along a curve, i. e. it is a reducing sphere for the original Heegaard splitting. □

Note that at the end of the proof above we have $\bar{S}$ dividing $M$ into two (not necessarily connected) 3-manifolds, $W_1$ and $W_2$. Each component of $W_i$ inherits a Heegaard splitting surface (a component of $S_i$) of lower genus than $S$. This splitting itself may be weakly reducible and we can continue the process. Ultimately an irreducible Heegaard splitting $M = H_1 \cup_S H_2$ is thereby broken up into a series of strongly irreducible splittings (see [ST3]). That is, we can begin with the handle structure determined by $H_1 \cup_S H_2$ and rearrange the order of the 1- and 2-handles, so that ultimately

$$M = M_0 \cup_{\bar{S}_1} M_1 \cup_{\bar{S}_2} \ldots \cup_{\bar{S}_m} M_m.$$

The 1- and 2-handles which occur in $M_i$ provide it with a strongly irreducible splitting (in each component) $A_i \cup_{P_i} B_i$ with $\partial_- A_i = \bar{S}_i$, $\partial_- B_{i-1} = \bar{S}_i$ for $1 \leq i \leq m$, $\partial_- A_0 = \partial_- H_1 \subset \partial M$, $\partial_- B_m = \partial_- H_2 \subset \partial M$. Each component of each $\bar{S}_i$ is a closed incompressible surface of positive genus and, for any $i$, only one component of $M_i$ is not a product. None of the compression bodies $A_i, B_{i-1}, 1 \leq i \leq m$ is trivial. If $\partial_- A$ or $\partial_- B$ is compressible in $M$ (so in particular $M$ is $\partial$-reducible) then respectively $A_0$ or $B_m$ is trivial (i. e. just a product). Such a rearrangement of handles will be called an *untelescoping* of the Heegaard splitting.

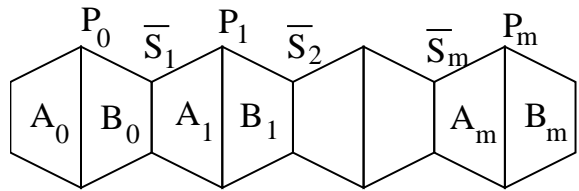

Figure 7.

So we see that just as a reducible splitting can be broken up by spheres into a connected sum of irreducible Heegaard splittings so irreducible Heegaard splittings can, by rearranging handles, be decomposed by incompressible surfaces into a sequence of strongly irreducible splittings. In effect, strongly irreducible splittings can be viewed as the fundamental building blocks of general Heegaard splittings.

The inverse process is also of interest. Suppose an incompressible surface $\bar{S}$ divides a connected 3-manifold $M$ into two pieces $M_0$ and $M_1$ and, for $i = 0, 1$ there are surfaces $P_i \subset M_i$ which divide (each component of) $M_i$ into compression bodies $A_i$ and $B_i$, with $\partial_- B_0 = \bar{S} = \partial_- A_1$. From this we can recover a Heegaard splitting of $M$ by a process called *amalgamation* (see [Sc1]). Informally, we regard the two Heegaard splittings as handle decompositions and rearrange the handles so that all the 1-handles are attached to $\partial_- A_0$ before the 2-handles are attached. More formally, do the following: The 3-manifold $B_0 \cup_{\bar{S}} A_1$ can be viewed as obtained from $\bar{S} \times [-1, 1]$ by attaching some 1-handles (from $B_0$) to $\bar{S} \times \{-1\}$ and some 1-handles (from $A_1$) to $\bar{S} \times \{1\}$. The attaching disks of these 1-handles, in $\bar{S} \times \{\pm 1\}$ can be taken to project to disjoint disks in $\bar{S}$. Collapse $\bar{S} \times I$ to $\bar{S}$. Then the 1-handles of $B_0$ are attached to $P_1 = \partial_+ B_1$, which makes it a compression body $B$, and the 1-handles of $A_1$ are attached to $P_0 = \partial_+ A_0$, which makes it a compression body $A$. Moreover, $\partial_+ A = \partial_+ B$. If we denote this surface $S$, then $M = A \cup_S B$ is a Heegaard splitting.

## 4. Heegaard splittings in nature: Seifert manifolds

An important class of irreducible 3-manifolds is the Seifert manifolds. We restrict our comments here to those Seifert manifolds constructed only with orientation preserving data, as it is these on which Heegaard splittings are best understood. We call these *fully orientable* Seifert manifolds.

Here is how a fully orientable Seifert manifold is constructed (see also [Sc]). Begin with $F$ a compact orientable surface (called the base surface of the Seifert manifold) and choose $n$ points $x_1, \ldots, x_n \in F$.

Choose small disjoint disks $E_i$ around the $x_i$. Let $F_- = F - \cup(\cup_i \partial E_i)$. In $F \times S^1$ do Dehn surgery on each $E_i \times S^1$ as follows (see [Boy] ): For each $1 \leq i \leq n$ remove $E_i \times S^1$ and glue back a solid torus $T_i$ so that a circle $\{x\} \times S^1 \subset \partial E_i \times S^1 \subset \partial F_i \times S^1$ is identified with a $(p_i, q_i)$ torus knot on $\partial T_i$. This means a knot going $p_i \geq 2$ times around the longitude (i. e. crossing a meridian $p_i$ times) and $q_i$ times around a meridian (i. e. crossing a longitude $q_i$ times). Because of ambiguity in the choice of longitude for $T_i$, $q_i$ is only defined mod $p_i$ and it is customary to take $1 \leq q_i < p_i$. Once this Dehn surgery is done, the projection of $\partial E_i \times S^1$ to $B^2$ extends to a projection of $T_i$ to $E_i$ in which the inverse image of each point of $E_i - \{x_i\}$ is a $(p_i, q_i)$ torus knot in $T_i$ and the inverse image of $x_i$ is the core circle of $T_i$. So the resultant manifold still projects to $F$. The inverse image of each point in $F$ is a circle (called a *fiber*). The inverse image of each $x_i$ is called an *exceptional fiber* and other fibers are called *regular fibers*.

The description of the Dehn surgery is not yet complete since there is still choice in how the cross-section $\partial E_i \times \{y\}$ is attached to $\partial T_i$. Any first choice could be altered by Dehn twists, in $\partial T_i$, along the $(p_i, q_i)$ torus knot image of $\{x\} \times S^1$, so the possible choices are parameterized by the integers. But there is also another ambiguity. There may be automorphisms of $F_- \times S^1$ which preserve the fiber structure (i. e. commute with projection to $F_-$). It is easy to see that, when $F$ has boundary, all of the choice of cross-section attachment can be absorbed into the ambiguity of what is the global cross-section $F_- \times \{y\}$, so the data above is sufficient to characterize the manifold up to homeomorphism that commutes with projection to $F$. If $F$ is closed, there is not so much flexibility of the choice of cross section so ultimately there is an integer's worth of choice involved in how the manifold is constructed. Typically this choice is realized by choosing another point $x_0 \in F_-$ and disk $E_0 \subset F_-$ containing it and doing Dehn surgery on $E_0 \times S^1$, gluing in a solid torus $T_0$ so that its longitude is identified with $\{x\} \times S^1 \subset \partial E_0 \times S^1$ and its meridian is identified with a cross-section of $\partial E_0 \times S^1$, of which there are an integer's worth of possibilities. The choice here determines what is called the *Euler number* of the Seifert manifold. (Working out the details in this paragraph is a good first step at understanding obstruction theory).

There are two ways in which the Seifert structure can induce natural Heegaard splittings on the 3-manifold and these are the subject of the next two sections. It is the principal result of [MS] that any irreducible splitting of a fully orientable Seifert manifold is one of these two types.

4.1. **Vertical splittings.** Suppose that $M$ is a fully orientable Seifert manifold, constructed as above, with base surface $F$, projection $p : M \to F$, and singular fibers the inverse images of $x_1, \ldots, x_n \in F$. Let $\Gamma$ be a connected graph in $F_-$ chosen so that

1. some nonempty subset of the $x_i, 1 \leq i \leq n$ are vertices of $\Gamma$
2. each component of $F - \Gamma$ is either a disk containing a single $x_i$ or an annulus containing a single boundary component. (But if $F$ is closed and $n = 1$, then $F - \Gamma$ is a disk not containing $x_1$.)
3. $\Gamma \cap \partial F$ consists of a (possibly empty) collection of boundary components $d_1, \ldots, d_d$.

Let $H_1 \subset M$ be the compression body whose spine is the union of $\{d_j \times S^1, j = 1 \ldots d\}$, a lift of $\Gamma$, and the singular fiber lying over each $x_i \subset \Gamma$. The complement of $H_1$ in $M$ is also a compression body, whose spine is the union of the boundary components of $M$ not in $\{d_j \times S^1\}$, the exceptional fibers not lying over $\Gamma$, and the lift of a "dual" complex to $\Gamma$. This creates a Heegaard splitting which is called a *vertical* Heegaard splitting.

It is not difficult to show that, up to isotopy, this construction is independent of $\Gamma$ but depends only on the choice of the $x_i$ that lie in $\Gamma$ and the choice of the boundary components $d_1, \ldots, d_d$.

**Theorem 4.1** ([Sc2]). *If $M$ is a fully orientable Seifert manifold and $\partial M \neq \emptyset$ then any irreducible Heegaard splitting is vertical.*

**Proof**: Here is a sketch of the complicated and ingenious argument. The proof is by induction on the number of exceptional fibers, with the case of no such fibers (i. e. $M = F \times S^1$) covered in [Sc1]. Let $e$ be an exceptional fiber, put in "thin" position with respect to a sweep-out (section 2.4) coming from the Heegaard splitting. This means, roughly, that if one considers how the circle fiber $e$ is intersected by the sweep-out, one cannot move the levels at which the maxima and minima of $e$ occur by pushing a maximum down and a minimum up until the the maximum is encountered before the minimum. This so simplifies $e$ that it can be moved to lie on a Heegaard surface in a way so that it intersects a meridian curve on one side in a single point. This is sufficient to ensure that $e$ can be made a core of a handle on one side, so that removing it from $M$ leaves the Heegaard surface as the splitting surface of $M - \eta(e)$. □

4.2. **Horizontal splittings.** Here is a specialized way to construct some fully orientable Seifert manifolds. Let $\hat{F}$ be an orientable compact surface and $h : \hat{F} \to \hat{F}$ be a periodic orientation preserving diffeomorphism, such that $h^n = identity$. Consider the mapping cylinder of $h$, a compact 3-manifold $M$ obtained by identifying, in $\hat{F} \times I$, each point $x \times 0 \in \hat{F} \times \{0\}$ with $h(x) \times 1$. Notice that $M$ is fibered by circles. For any point $x \in \hat{F}$ the union of the images of $\{h^i(x)\} \times I \subset \hat{F} \times I, 1 \leq i \leq n$ is a circle which typically intersects a cross-section $\hat{F} \times \{s\}$ in $n$ points. For some discrete (hence finite) set of points in $\hat{F}$, the orbit may be of length only a proper factor $l$ of $n$ and then the corresponding circle intersects a cross-section only in $l$ points. It's easy to see that this gives $M$ a Seifert manifold structure in which the base surface is $F = \hat{F}/h$, a surface over which $\hat{F}$ is a branched covering.

A Seifert manifold can be given such a structure if and only if its Euler number is zero (see [Sc]). In particular, any fully orientable Seifert manifold with non-empty boundary can be given such a structure. Suppose $M_-$ is a fully orientable Seifert manifold whose boundary is a single torus $T$. View $M_-$ as the mapping cylinder of a diffeomorphism $h : \hat{F} \to \hat{F}$, where $\hat{F}$ is a compact orientable surface with a single boundary component. Now $\partial \hat{F} \subset \partial M_-$ is a circle $c$ transverse to the fibers of $M_-$. Attach a solid torus $T$ to $\partial M_-$ so that a longitude goes to $c$. (There is an integer's worth of choice of how the meridian is attached.) This creates a closed Seifert manifold $M$ which can be split into two pieces: The image of $\hat{F} \times [0, 1/2]$ and the union of $T$ and $\hat{F} \times [1/2, 1]$. Since both pieces are homeomorphic to $\hat{F} \times I$ (in the latter case because $T$ becomes just a collar of $\partial \hat{F} \times I$), each is a handlebody. Thus we get a Heegaard splitting, and this construction is called a *horizontal* splitting of $M$.

**Theorem 4.2** ([MS],[Sc3]). *An irreducible Heegaard splitting of a fully orientable Seifert manifold is either horizontal or vertical.*

**Proof**: A special argument ([Sc3]) is needed for small Seifert manifolds; we sketch here the proof when $M$ contains an essential vertical torus. Suppose $H_1 \cup_S H_2$ is the irreducible splitting and suppose that it is weakly reducible. Since the splitting is irreducible it follows from the proof of Theorem 3.11 that, if $S$ is maximally and independently compressed in both directions, the result is an incompressible surface $\bar{S}$. Furthermore $S$ can be viewed as assembled from Heegaard splittings of $M - \bar{S}$ by amalgamating along $\bar{S}$. Any incompressible surface in a Seifert manifold can be isotoped to be either a collection of vertical tori, or the fiber in a fibering of $M$ over $S^1$.

Suppose $\bar{S}$ is a collection of vertical tori. Then $M - \bar{S}$ is a Seifert manifold with boundary, and so any Heegaard splitting is vertical, by Theorem 4.1. Thus $S$ is obtained by amalgamating vertical splittings along vertical tori, from which it follows immediately that $S$ is also vertical.

Suppose $\bar{S}$ is a set of fibers of a fibering over $S^1$. Then it splits $M$ into pieces of the form $\bar{S} \times I$, and induces a Heegaard splitting on each piece. Heegaard splittings of $surface \times I$ are well-understood ([ST1]) and examination shows that in fact the compressing could have been done in such a way that $\bar{S}$ would be a collection of vertical tori, reducing to the previous case.

Suppose finally that $S$ is strongly irreducible. An argument similar in spirit to the thin position argument of Theorem 4.1 proves that a fiber $f$ can be isotoped onto the surface $S$. Then the Seifert manifold $M_- = M - \eta(f)$ is split in two pieces by the surface $S_- = S - \eta(f)$. If $S_-$ is incompressible in $M_-$ then it is the fiber of a fibration of $S_-$ over $S^1$ and it follows easily that the original $S$ is a horizontal splitting. If $S_-$ is compressible then, since $S$ is strongly incompressible, after a maximal number of compressions into one handlebody, say $H_2$, $S_-$ becomes an incompressible surface $S*$ in $M_-$.

If $S^*$ is a vertical annulus then the union of $S^*$ and $S \cap \eta(f)$ is a torus in $H_1$ so it bounds a solid torus. The core of the torus is a fiber and the manifold obtained by deleting it has $S$ as a Heegaard splitting surface. It follows from 4.1 that $S$ is vertical.

If $S^*$ is the fiber of a fibering of $M_-$ over $S^1$ then $S^*$ splits $M_-$ into handlebodies and $S$ is a further Heegaard splitting of one them. But any (non-trivial) splitting of a handlebody is stabilized, hence reducible (from Proposition 3.6 and Theorem 3.7) and $S$ is assumed irreducible. □

On the other hand, not all vertical and horizontal splittings are irreducible. Exactly which ones are has been worked out in [Se2].

## 5. CONNECTIONS WITH GROUP PRESENTATIONS

In this section, assume that $M = H_1 \cup_S H_2$ is a closed manifold, and hence that both $H_1$ and $H_2$ are handlebodies, say of genus $g$. This implies that $\pi_1(H_1)$ is a free group on $g$ generators. A choice of basepoint and a complete collection $\Delta = \{D_1, \ldots, D_g\}$ of oriented meridian disks determines a presentation of $\pi_1(H_1)$, namely, for any based loop in $H_1$, write down $x_i$ every time the loop passes through the disk $D_i$ in a direction consistent with its normal orientation and $x_i^{-1}$ if the direction is inconsistent. Similarly, a complete collection $E_1, \ldots, E_s, s \geq g$ of meridian disks for $H_2$ then determines a presentation of $\pi_1(M)$. Each curve $\partial E_k$, when viewed as a (conjugacy class) in $\pi_1(H_1)$, and so as a word $r_k$ in $\{x_i\}$, is a relator for the fundamental group. That is, $\pi_1(M)$ has the presentation $\{x_1, \ldots, x_g; r_1, \ldots, r_s\}$. We say that this presentation is *geometrically realized*.

How much does this presentation depend on choices made? We'll restrict attention to $H_1$ (which yields the generators) since the situation is much the same for the relators.

**Lemma 5.1.** *Any set $\{y_1, \ldots, y_g\}$ of generators of $\pi_1(H_1)$ can be geometrically realized.*

**Proof**: There are a specific set of moves on generators, called the Nielsen moves, which will transform a given geometrically realized set of generators $\{x_1, \ldots, x_g\}$ into $\{y_1, \ldots, y_g\}$. But an examination of these moves (see e. g. [MKS, 3.1] ) shows that each move can be realized by a geometric move, either sliding one 2-handle over the other (i. e. band-summing one meridian disk to another) or reversing the orientation of a disk, or just naming the disks in a different order. □

Motivated in part by this lemma it makes sense to introduce the following definition.

**Definition 5.2** ([LM]). *Two generating systems $\{u_1, \ldots, u_g\}$ and $\{v_1, \ldots, v_g\}$ for $\pi_1(M)$ are* Nielsen equivalent *if there is an epimorphism $\phi : F_g \to \pi_1(M)$ and bases $\{x_1, \ldots, x_g\}$ and $\{y_1, \ldots, y_g\}$ for the free group $F_g$, such that $\phi(x_i) = u_i$ and $\phi(y_i) = v_i$.*

Less formally, if we view a presentation of $G$ as an epimorphism of the free group (with specified generators) onto $G$, then two presentations are Nielsen equivalent if there is an automorphism of the free group which realizes the change in specified generators.

It follows that a Heegaard splitting $M = H_1 \cup_S H_2$ (and a choice of which handlebody is $H_1$) specifies a single Nielsen equivalence class of presentations of $\pi_1(M)$. For if $\{u_1, \ldots, u_g\}$ and $\{v_1, \ldots, v_g\}$ are the generating systems induced by different choices of meridian disks for $H_1$, then in the definition above substitute $\pi_1(H_1)$ for $F_g$, let the inclusion induce $\phi$, and deduce that the presentations are Nielsen equivalent.

To generalize slightly:

**Theorem 5.3.** *If two Heegaard splittings $H_1 \cup_S H_2$ and $H_1' \cup_{S'} H_2'$ of the same closed $3$-manifold $M$ are isotopic then (for the appropriate choice of $H_1$ and $H_1'$) their corresponding geometrically realized presentations are Nielsen equivalent.*

**Proof**: Inner automorphism is a Nielsen equivalence. □

This gives a powerful algebraic tool to show that Heegaard splittings are not isotopic.

The structures of Heegaard splittings discussed in section 3 above have implications for the induced group presentations. For example,

if a Heegaard splitting is stabilized, then there is a Nielsen equivalent presentation in which a relator is precisely a generator. If it is reducible, then there is a Nielsen equivalent presentation which splits as a free product of two presentations. (The genus one splitting of $S^1 \times S^2$ and also the non-orientable 2-sphere bundle over $S^1$ are exceptions. And of course the groups presented might be trivial.) If it is weakly reducible, then there is a Nielsen equivalent presentation in which at least one generator does not appear in at least one relator.

Not all presentations can be geometrically realized. For example, Boileau and Zieschang [BZ] have shown that certain Seifert manifolds have fundamental groups which admit presentations with two generators, whereas the minimal genus of any Heegaard splitting (hence the rank of any geometrically realized presentation) is at least three. Montesinos [Mn] has used this example to show that a presentation Nielsen equivalent to a geometric presentation may not be geometric.

For details on this and other examples, see [Zi].

## 6. Uniqueness

How many distinct Heegaard splittings does a 3-manifold have? We have already seen that any Heegaard splitting can be stabilized, so the question only becomes interesting if we restrict to Heegaard splittings which are not stabilized.

6.1. **The 3-sphere.** In 1968 Waldhausen [Wa] showed that any positive genus Heegaard splitting of $S^3$ is stabilized (3.7), so that the only genus $g$ splitting is the obvious one, obtained by stabilizing $q$ times the splitting of $S^3$ into 3-balls. This was the first uniqueness result. Here is a sketch of a later proof ([ST2], [Ot]).

**Theorem 3.7.** *Any positive genus Heegaard splitting of $S^3$ is stabilized.*

**Proof**: Suppose $S^3 = H_1 \cup_S H_2$ and $\Sigma$ is a spine of $H_1$. We may assume $\Sigma$ is a tri-valent graph in $S^3$ and we are allowed to do edge-slides. Choose a Morse function $h : S^3 \to [-1, 1]$ which has a single minimum (at height -1) and a single maximum (at height 1) and which restricts to a Morse function on $\Sigma$. Put $\Sigma$ in "thin position" with respect to this height function. In outline, this means that you can't push down a maximum (this includes trivalent vertices in which two edges leave the vertex from below) so that it moves below a minimum (this includes trivalent vertices in which two edges leave the vertex from above) without introducing new critical points.

It suffices to show there is an unknotted cycle $\gamma \subset \Sigma$. For then $S$ would also be a Heegaard splitting surface for the solid torus $S^3 - \eta(\gamma)$.

This splitting would necessarily be boundary reducible (Proposition 3.6) which means that the original splitting $S$ was stabilized.

Consider a collection $\Delta \subset S^3$ of meridian disks of $H_2$, extended into $H_1$, so that its interior is embedded in $S^3 - \Sigma$ and its (singular) boundary lies in $\Sigma$. The first observation is that we may as well assume $\partial \Delta$ runs across every edge of $\Sigma$, for otherwise $H_1 \cup_S H_2$ would be reducible (Theorem 3.11). If the splitting were reducible then a reducing sphere splits $S$ into two Heegaard splittings of $S^3$ each of smaller positive genus, and we would be done by induction.

Consider when a level sphere $S_t = h^{-1}(t)$ cuts off from $\Delta$ a subdisk sufficiently simple that it can be used to slide part of an edge of $\Sigma$ so that it lies on $S_t$. It's easy to see that this is true just below the highest point of $\Sigma$ and just above the lowest point. In the former case the disk can be used to lower the maximum slightly and in the latter to raise the minimum. Suppose we simultaneously (i. e. for the same level sphere) have two subdisks of $\Delta$, one of which lowers a maximum and the other of which raises a minimum. Then either this violates thin position (when we can push the maximum slightly lower without interfering with the minimum) or the two edges which we have pushed onto the level sphere have the same ends, i. e. they create an unknotted cycle and we are done.

We know then that a sufficiently high sphere cuts off a subdisk of $\Delta$ lowering a maximum, a sufficiently low sphere cuts off a subdisk raising a minimum and, if subdisks of both types are cut off simultaneously, then we are done. So it suffices to eliminate the possibility that neither type occurs, that is, there is a height $t_0$ so that no subdisk cut off by $S_{t_0}$ from $\Delta$ can be used either to raise a minimum or lower a maximum. But this situation cannot in fact occur, by an argument reminiscent of the second proof of Theorem 3.4, with $S_{t_0}$ playing the role of the reducing sphere. □

6.2. **Seifert manifolds.** One might have hoped that this situation would generalize - that any compact 3-manifold would have ( up to stabilization) a unique Heegaard splitting. In 1970 R. Engmann ([En], see also [Bi]) showed that the connected sum of certain pairs of Lens spaces could have two non-homeomorphic Heegaard splittings of genus two (hence not stabilized). Examples were shortly found of prime manifolds with the same property ([BGM]). A rather spectacular generalization by Lustig and Moriah is the main theorem of [LM]. It is a good illustration of the usefulness of Theorem 5.3, so we sketch the central idea.

Let $M$ be a fully orientable Seifert manifold, constructed as in section 4, with base surface $F$, projection $p : M \to F$, and singular fibers the inverse images of $x_1, \ldots, x_n \in F$. Given details of the fibering around the $x_i$ and the Euler number of $M$ it is straightforward to write down a presentation of $\pi_1(M)$. It's easy to see directly that the element $h \in \pi_1(M)$ represented by a regular fiber is central in $\pi_1(M)$.

Consider the quotient group $G = \pi_1(M)/<h>$. The complement $M_-$ of the exceptional fibers is $F_- \times S^1$ so the effect of factoring out $<h>$ is to reduce $\pi_1(M_-)$ to $\pi_1(F_-)$. In the solid torus surrounding an exceptional fiber a meridian crosses a fiber some $p_i \geq 2$ times. The effect is to kill the $p_i$ multiple of $\partial E_i$ in $\pi_1(F_-)$. The upshot is that $G$ is a Fuchsian group and in particular has a faithful presentation into $PSL_2(\mathbf{C})$. This special structure provides an extra tool for determining when group presentations are Nielsen equivalent. (Note that Nielsen equivalent presentations of $\pi_1(M)$ descend to Nielsen equivalent presentations of $G$.)

This extra information is sufficient to show that, in most cases, two vertical splittings of the same fully orientable 3-manifold are isotopic only if the equivalence is more or less obvious, e. g. the invariants of the exceptional fibers lying in the graph $\Gamma$ are the same (see section 4.1). In particular this leads to a complete classification of irreducible Heegaard splittings of most fully orientable Seifert 3-manifolds with boundary (see [Sc2]). In the case of closed Seifert manifolds, there is still some puzzlement about how horizontal splittings fit into the classification scheme. For example, whereas a vertical splitting of a closed Seifert manifold with base surface of genus $g$ and with $k$ exceptional fibers is $2g + k - 1$, there are some such Seifert manifolds (over $S^2$, i. e. $g = 0$) which have horizontal splittings of genus $k - 2$.

### 6.3. Genus and the Casson-Gordon examples.

The last comment prompts the following question: Do we at least know that all irreducible splittings of the same 3-manifold have the same genus? In 1986 Casson and Gordon gave an example which shows that the answer is an emphatic no [CG2], [Ko]. What they show is that there is a closed orientable 3-manifold (in fact infinitely many) which has irreducible splittings of arbitrarily high genus. We outline the construction.

Begin with the following fact [Pa] : There are certain pretzel knots $k \subset S^3$ with the property that they have incompressible Seifert surfaces of arbitrarily high genus (these are explicitly constructed) and for each of these surfaces the complement in $S^3$ is a handlebody. Pick one of these knots, and let $F_n$ be an incompressible Seifert surface of genus $n$ whose complement in $S^3$ is a handlebody. Then $S^3 = \eta(F_n) \cup (S^3 -$ $int(\eta(F_n)))$ is a (highly reducible) genus $2n$ Heegaard splitting of $S^3$, and $k$ is isotopic to a curve on the splitting surface $S = \partial \eta(F_n)$. Let $M_q$ be the 3-manifold obtained by doing $1/q$ surgery on $k$ ($q$ an integer). One way to view this is to imagine pulling the two handlebodies apart along a strip parallel to $k \subset S$ then gluing the two strips back together via a $q$-fold Dehn twist. So in particular the construction naturally gives a genus $2n$ Heegaard splitting of $M_q$. For $q$ a large integer ($q \geq 6$ suffices) it turns out (see below) that the resulting splitting is strongly irreducible. Thus a specific $M_q$ will have splittings, built as above for different values of $n$, of arbitrarily high genus.

The critical ingredient in the above argument is then

**Theorem 6.1** ([CG2]). *Suppose $M = H_1 \cup_S H_2$ is a weakly reducible Heegaard splitting of the closed manifold $M$. Let $k$ be a simple closed curve on the splitting surface $S$ so that $S - \eta(k)$ is incompressible in both $H_i, i = 1, 2$. Let $M_q$ be the manifold obtained by $1/q$ surgery on $k$. Then for $q \geq 6$ the associated Heegaard splitting (induced as above) on $M_q$ is strongly irreducible.*

See [MS, Appendix] for a proof. The idea is this: $k$ necessarily intersects any meridian disk on either side, since $S - \eta(k)$ is incompressible on both sides. Sufficient Dehn twisting along $k$ then will stretch any meridian of one side so that it intersects any meridian disk of the other.

In fact ([Ko]) the number of Heegaard splittings at each even genus is bounded below by a polynomial in the genus.

6.4. **Other uniqueness results.** We briefly note that there are other manifolds which are known to have unique irreducible Heegaard splittings (for a particular distribution of boundary components between $H_1$ and $H_2$). A perhaps not exhaustive list is the following:
- $S^2 \times S^1$ [Wa]
- Any Lens space [Bo] [BoO] (see also 7.10)
- Any *(closed orientable surface)* $\times I$ [ST1] [BO1]
- Any *(compact orientable surface)* $\times S^1$ [Sc1] [BO1]

## 7. The stabilization problem

We noted in section 2.1 that every compact 3-manifold admits a Heegaard splitting, since every 3-manifold has a triangulation. Similarly, since any two triangulations of the same 3-manifold are PL equivalent (see [Mo], [Bn]) it follows that any two Heegaard splittings have a common stabilization. The classical argument, which goes back to Reidemeister and Singer, is more complicated than one might expect. See [AM] for details. A more straightforward proof has recently been noted by Fengchun Lei ([L]). It exploits the fact that the new proofs of the uniqueness of splittings of $S^3$ (see Theorem 3.7) do not require the Reidemeister-Singer result (as Waldhausen's original proof did).

**Theorem 7.1** ([L]). *Any two Heegaard splittings of the same compact 3-manifold have a common stabilization.*

**Proof**: The case in which $M$ is closed is representative. First note that Waldhausen's theorem (Theorem 3.7) combined with Theorem 3.6 easily show (by induction on genus) that any Heegaard splitting of a handlebody is either trivial or stabilized. So then suppose $A \cup_P B$ and $X \cup_Q Y$ are two splittings of the closed manifold $M$. We can assume that the spines $\Sigma_A$ and $\Sigma_Y$ of $A$ and $Y$ are disjoint in $M$; let $W$ be the compact manifold obtained by removing an open neighborhood $\eta(\Sigma_A \cup \Sigma_Y)$ and let $S$ be a Heegaard splitting surface for $W$. Then $S$ is also a Heegaard splitting surface for $M$. Furthermore, $S$ is also a splitting surface for the handlebody $B = M - \eta(\Sigma_A)$ so it stabilizes $A \cup_P B$ and for the handlebody $X = M - \eta(\Sigma_Y)$ so it stabilizes $X \cup_Q Y$. □

We noted in section 6.2 that the connection between Heegaard splittings and group presentations and the known structure of Heegaard splittings of Seifert manifolds show that some manifolds have distinct irreducible Heegaard splittings. This raises the natural question: How much do we need to stabilize before two splittings of the same 3-manifold become isotopic? More generally, now that we know that a manifold can have quite different Heegaard splittings, how can such distinct Heegaard splittings be compared?

As a cautionary tale, revealing the depth of our ignorance on the first question, consider the large gap between what is known and what is not:

**Theorem 7.2** ([Sc3]). *Two irreducible Heegaard splittings of the same fully orientable Seifert manifold have a common stabilization requiring, for one of the splittings, at most one stabilization.*

**Theorem 7.3** ([Se1]). *For the irreducible Heegaard splittings of $M_q$ constructed in section 6.3, any two Heegaard splittings of the same $M_q$ have a common stabilization requiring, for one of the splittings, at most one stabilization.*

In fact, there is no example of distinct Heegaard splittings of the same closed 3-manifold which cannot be made isotopic by a single stabilization of one of the splittings, and sufficient stabilizations of the other to ensure that the genus of the two surfaces is the same. One could thus make the very optimistic

**Conjecture 7.4.** *Suppose $H_1 \cup_S H_2$ and $H'_1 \cup_{S'} H'_2$ are Heegaard splittings of the same 3-manifold of, genus $g \leq g'$ respectively. Then the splittings obtained by one stabilization of $S'$ and $g' - g + 1$ stabilizations of $S$ are isotopic.*

At the other extreme are two theorems which put limits on how much stabilization is needed, in terms of the genera of the two original splittings.

**Theorem 7.5** ([Jo], Theorem 31.9). *Suppose $M$ is a Haken 3-manifold containing no non-trivial essential Stallings fibrations. Then the number of stabilizations required to guarantee that a genus $g$ splitting of $M$ is isotopic to a genus $g'$ splitting is some polynomial function (perhaps depending on $M$) of $g$ and $g'$.*

The gap here is rather huge. An ideal sort of theorem would be one which gives an explicit bound, independent of the manifold, on the number of stabilizations required, expressed in terms of the genera of the splittings being considered. This would be an important step toward solving the "homeomorphism problem" - find an algorithm which will determine if two compact 3-manifolds are homeomorphic - because it would reduce the problem to the case in which there are isotopic Heegaard splittings of the same known genus for the two manifolds.

In this direction is the following theorem, which applies to all irreducible splittings of compact orientable non-Haken 3-manifolds:

**Theorem 7.6** ( [RS1], Theorem 11.5 ). *Suppose $X \cup_Q Y$ and $A \cup_P B$ are strongly irreducible Heegaard splittings of the same closed orientable 3-manifold $M$ and are of genus $p \leq q$ respectively. Then there is a genus $8q + 5p - 9$ Heegaard splitting of $M$ which stabilizes both $A \cup_P B$ and $X \cup_Q Y$.*

It appears that a similar explicit, but quadratic, bound can be found for Haken 3-manifolds using [RS3] and [RS2]. The former extends Theorem 7.6 to the bounded case. The machinery of the latter, [RS2], illustrates how to extend certain general position arguments of [RS1] to weakly reducible splittings that have been untelescoped, as described in section 3.9.

The proof of Theorem 7.6 is quite complicated, but a crucial ingredient is a theorem that describes how two strongly irreducible splittings can be moved to intersect in a way that contains much information about both splittings.

**Theorem 7.7** ( [RS1], Theorem 6.2 ). *Suppose $X \cup_Q Y$ and $A \cup_P B$ are strongly irreducible Heegaard splittings of the same closed orientable* 3-manifold $M \neq S^3$. Then $P$ and $Q$ can be isotoped so that $P \cap Q$ is a non-empty collection of curves which are essential in both $P$ and $Q$.

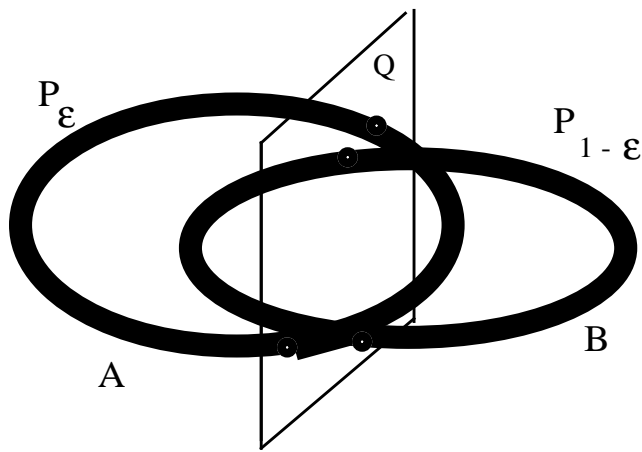

Figure 8.

The proof is an application of Heegaard splittings as sweep-outs (section 2.4). For $A \cup_P B$ a Heegaard splitting of a closed 3-manifold $M$, let $\Sigma_A$ and $\Sigma_B$ be spines of the handlebodies $A$ and $B$ respectively. Recall

**Definition 7.8.** *A* sweep-out *associated to the Heegaard splitting $A \cup_P B$ is a relative homeomorphism $H : P \times (I, \partial I) \to (M, \Sigma_A \cup \Sigma_B)$ which, near $P \times \partial I$, gives a mapping cylinder structure to a neighborhood of $\Sigma_A \cup \Sigma_B$.*

*Given such a sweep-out $H$ and a value $s, 0 \leq s \leq 1$, let $P_s$ denote $H(P \times s)$, $P_{<s}$ denote the handlebody $H(P \times [0, s])$ and $P_{>s}$ denote the handlebody $H(P \times [s, 1])$. Note that $P(s), s \neq 0, 1$ is a copy of the splitting surface, $P_0 = \Sigma_A$ and $P_1 = \Sigma_B$.*

Consider how the surfaces $P_s$ intersect a distinct Heegaard splitting surface $Q$ in $M = X \cup_Q Y$. Assume $Q$ is in general position with respect to $\Sigma_A \cup \Sigma_B$ (so the spines intersect $Q$ transversally in a finite number of points) and the sweep-out $H$ is generic with respect to $Q$. Then, for small values of $\epsilon$, $P_{<\epsilon}$ is very near $\Sigma_A$, so $P_{<\epsilon} \cap Q$ is a (possibly empty) collection of meridian disks of $A$. Symmetrically $P_{>1-\epsilon}$ is very near $\Sigma_B$, so $P_{>1-\epsilon} \cap Q$ is a (possibly empty) collection of meridian disks of $B$. Throughout the sweep-out, at least generically, $P_s \cap Q$ is a disjoint collection of simple closed curves in $Q$. (See fig. 8.)

Note that $P_s \cap Q$ cuts off in $Q$ meridian disks for $A$ when $s$ is small, meridian disks for $B$ when $s$ is large and can't cut off simultaneously meridian disks for both, since $A \cup_P B$ is strongly irreducible. It follows that for some value of $s$, no meridian is cut off. That is (with a minor amount of additional fuss) every curve of $P_s \cap Q$ is essential in $Q$.

In order to prove Theorem 7.7 we would like to apply a similar argument simultaneously to sweepouts $P_s, Q_t$ of $M$ corresponding to the different Heegaard splittings of $M$. Cerf theory (see [C]) can be used to make the following informal remarks rigorous. A good way to think visually of the discussion below is to consider the surfaces $P_s$ and $Q_t$ as parameterized by the point $(s,t)$ in the square $I \times I = \{(s,t) | 0 \leq s, t \leq 1\}$.

Away from $\partial(I \times I)$, four things can happen:

- At a generic value of $(s,t)$, $P_s$ and $Q_t$ intersect transversally in a collection of simple closed curves $c_{(s,t)}$ which we can regard as lying in either $P \cong P_s$ or $Q \cong Q_t$.

- On a one-dimensional stratum of $I \times I$, $P_s$ and $Q_t$ intersect transversally except at a single non-degenerate tangency point. A good way to think about this is to begin at a value of $(s,t)$ at which there is such a tangency point. Now imagine letting $s$ ascend (or descend) at just the rate required to ensure that the tangency point persists as $t$ ascends. This requirement defines $s$ as a function of $t$, and so parameterizes an arc inside the square $I \times I$. (Note that the slope of the arc is positive or negative depending on whether the ascending normal vectors to $P_s$ and $Q_t$ are parallel or anti-parallel. Thus the sign of the slope is fixed, providing a surprising order to the picture. So far, this additional order has not proven useful.)

- A discrete set of points $(s,t)$ for which $P_s$ and $Q_t$ have exactly two non-degenerate points of tangency but are otherwise transverse. For example, as $(s,t)$ traces out the arc as just described, there may be points of tangency which occur elsewhere. These are the discrete critical points of double tangency.

- A discrete set of points at which $P_s$ and $Q_t$ intersect transversally except for a single degenerate tangent point (locally modelled on $P_s = \{(x,y,z) | z = 0\}$ and $Q_t = \{(x,y,z) | z = x^2 + y^3\}$). These are so-called "birth-death" points, and play no important role in our discussion.

The set of points $(s,t)$ at which the intersection is non-generic forms a 1-complex $\Gamma$ called the *graphic* of the sweep-outs in the interior of $I \times I$. The graphic $\Gamma$ naturally extends to a properly imbedded 1-complex in all of $I \times I$: A point $(0,t)$, say, on $\{0\} \times I \subset \partial(I \times I)$ represents simultaneously the spine $\Sigma_A$ of handlebody $A$ (since $s = 0$) and the surface $Q_t$. Generically these are transverse, implying that $P_\epsilon$ and $Q_t$ are transverse for $\epsilon$ small. There are two types of exceptions: For finitely many values of $t$, $\Sigma_A$ is tangent to $Q_t$ at a single point in the interior of one of its edges. At finitely many other values of $t$, $Q_t$ crosses a vertex of $\Sigma_A$. It's easy to see that, in both these cases, nearby interior points are in the graphic, and vice versa, so $\Gamma$ extends to a graphic in the closed square. (See fig. 9.)

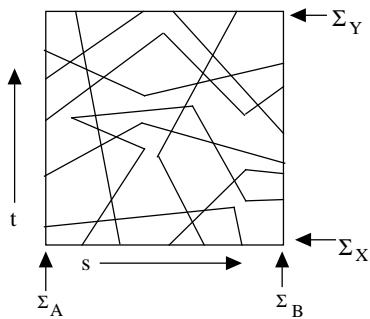

Figure 9.

**Proof of Theorem 7.6:** The graphic $\Gamma$ cuts $I \times I$ up into regions, in each of which the curves $c_{(s,t)}$ vary only by an isotopy in $P$ and $Q$. In some regions (e. g. near $\{0\} \times I$), $c_{(s,t)}$ cuts off meridian disks for $A$ lying in $Q$. In other regions (e. g. near $\{1\} \times I$), $c_{(s,t)}$ cuts off meridian disks for $B$ lying in $Q$. One can't have both occur in the same region, since $P$ is strongly irreducible. (Indeed they can't even occur in adjacent regions, but this is not immediately obvious.) Similarly there are regions in which $c_{(s,t)}$ cuts off a meridian of $X$ or $Y$, but not both, in $P$. We now apply a "mountain-pass" sort of argument: Given what we have described, there must be some point in the interior of $I \times I$ in which $c_{(s,t)}$ cuts off no meridians whatsoever. Such a point is the point we seek. It corresponds to an intersection in which $c_{(s,t)}$ is non-empty and each curve is essential in both $P$ and $Q$. (See fig. 10.)

What is not apparent in the above argument is why the point $(s,t)$ we have located is a generic point, nor is it clear how we can guarantee that the intersection $c_{(s,t)}$ is non-empty at this point. The details here require close argument, see [RS1]. □

To illustrate the power of this argument, we classify the irreducible splittings of the Lens space ([Bo], [BoO]). But first observe that (with a bit of reorganization) the argument above shows that any positive genus Heegaard splitting of $S^3$ is stabilized: Compare such a splitting with the index zero splitting (i. e. by $S^2$) of $S^3$. How could a curve of intersection with $S^2$ be essential in $S^2$?

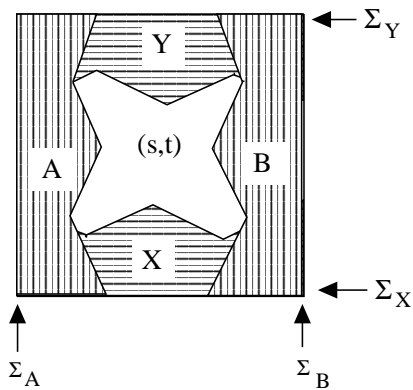

FIGURE 10.

**Corollary 7.9** ([Bo]). *Any two genus one Heegaard surfaces in a lens space are isotopic.*

**Proof**: Let $P$ and $Q$ be two genus one Heegaard surfaces in a lens space, separating the lens space, as usual, into solid tori $A$ and $B$ and solid tori $X$ and $Y$ respectively. $P$ and $Q$ may be isotoped so that they intersect in a non-empty family of essential circles, at which point $c_{(s,t)}$ cuts each up into annuli. One can pass annuli of $P$ parallel to annuli of $Q$ through each other until only two curves of intersection remain. At this point it is easy to show that the remaining annuli of $P - Q$ are parallel to the annuli $Q - P$. This means $P$ is parallel to $Q$. □

**Corollary 7.10** ([BoO]). *Any irreducible Heegaard splitting of a lens space has genus one.*

**Proof**: Let $A \cup_P B$ be a genus one Heegaard splitting of a lens space $L$ and $X \cup_Q Y$ be a splitting of higher genus. Since $L$ contains no incompressible surfaces, it suffices to show that $Q$ is weakly reducible. $P$ and $Q$ may be isotoped so that they intersect in a non-empty family of essential circles. There are two cases:

   **Case 1**: $Q \cap A$ and $Q \cap B$ both contain components which aren't annuli.

   As above, remove parallel annuli in $P$ and $Q$ by an isotopy. Then one can show, by $\partial$-compressing $Q$ in both $A$ and $B$, that somewhere in $A$ there is a meridian of $X$ and somewhere in $B$ a meridian of $Y$ (or vice versa). This contradict the strong irreducibility of $Q$. So we are reduced to **Case 2**: $Q \cap A$ or $Q \cap B$ (say the former) consists entirely of annuli.

When we remove parallel annuli, one can show that in the end, $Q$ actually lies in the solid torus $B$ and in fact induces a Heegaard splitting of $B$. But it follows from Proposition 3.6 that any higher genus Heegaard splitting of a solid torus is stabilized. □

## 8. Normal surfaces and decision problems

### 8.1. Normal surfaces and Heegaard splittings.

Much of our understanding of 3-manifolds depends on the surfaces they contain. Their most elementary taxonomy is expressed by reference to these surfaces: $M$ is irreducible if it contains no essential sphere, it's Haken if it contains a higher genus incompressible surface, it's atoroidal (and so, if closed and Haken, hyperbolic) if no such incompressible surface is a torus. Heegaard splittings have proven fundamental both in understanding the behavior of these surfaces and in developing algorithms (typically impractical) for classifying 3-manifolds within this taxonomy.

We briefly review the theory of normal surfaces [Ha2]. A good source for more detail is [JR].

Let $M = H_1 \cup_S H_2$ be a Heegaard splitting of a closed manifold (the case where $M$ is merely compact is an easy variation). Regard each $H_i$ as a handlebody, the union of some 0-handles (one would suffice) and some 1-handles. Suppose $F$ is a closed surface in $M$. Then $F$ can be isotoped so that it is disjoint from the points which are the cores of the 0-handles of $H_2$ and intersects transversally each of the cores of the 1-handles. By thickening these cores to the full handle-body we can isotope $F$ so that it intersects $H_2$ in some finite number of copies of 2-disk cocores $\Delta_2$ (meridians) of the 1-handles of $H_2$. Call the number of such disks in $F$ the *weight* of $F$.

Consider how $F$ then intersects $H_1$. It is helpful to recall the discussion of Heegaard diagrams in section 2.3. The handlebody $H_1$, when cut up by a family $\Delta_1$ of meridians, becomes a collection of 3-balls, the 0-handles of $H_i$. With little loss of generality we will assume in this discussion that there is a single 3-ball, $B^3$. In $\partial B^3$, the attaching curves of $\Delta_2$ become a 1-manifold $\mathcal{A}$ in $\partial B^3 - V$. The arcs are regarded as edges of a graph $\Gamma$ whose vertices $V$ correspond two-to-one to the meridians $\Delta_1$ of $H_1$. (We will here expand $\Gamma$ to include the simple closed curves of $\mathcal{A}$.) We may as well assume that $F$ is transverse to $\partial B^3$, so that $F \cap B^3$ is a properly imbedded surface lying in $B^3$. Because we have already assured that $F \cap H_2$ consists of copies of $\Delta_2$ we know that the collection of simple closed curves $F \cap \partial B^3$ is the union of parallel copies of components of $\mathcal{A}$ outside of $V$ and, inside of $V$, consists of some properly imbedded 1-manifold.

**Definition 8.1.** *The surface $F$ is* normal *with respect to $H_1 \cup_S H_2$ if*
1. *Each component of $F \cap B^3$ is a disk.*
2. *No component of $F \cap \partial B^3$ lies entirely in a fat vertex $v$.*
3. *Each component of $F \cap \partial B^3$ contains at most one copy of any edge of $\Gamma$.*

**Definition 8.2.** *A property of surfaces in 3-manifolds is called* compression preserved *if whenever a surface $F$ in $M$ has this property, and $F'$ is obtained from $F$ by a 2-surgery, then some set of components of $F'$ (not inessential spheres) also has this property.*

Examples of compression preserved properties are
- $F$ is a reducing sphere.
- $F$ is an injective surface (i. e. $\pi_1(F) \to \pi_1(M)$ is injective).
- $F$ has maximal Euler characteristic in its homology class (ignoring inessential spheres).

We then have:

**Theorem 8.3.** *If a closed 3-manifold $M$ contains a surface with a compression preserved property, then it contains a normal surface with the same compression preserved property.*

**Proof**: Choose a surface in $M$ which has the property and also has minimal weight. By compressing along disks lying in $V$ we can remove any components of $F \cap \partial B^3$ that lie completely in $V$. By compressing along disks lying slightly inside $\partial B^3$ we can arrange that the surface intersects $B^3$ in disks and in components lying entirely inside $B^3$. The latter can be discarded since they compress to inessential spheres and so, by Definition 8.2, they do not have the property. Finally, if any component of $F \cap \partial B^3$ contains more than one arc parallel to an edge $\gamma$ in $\Gamma$ then there is a $\partial$-compression of the corresponding disk in $F \cap B^3$ to an arc in $\eta(\gamma) \subset \partial B^3$. Then push across a 2-handle, reducing the weight of $F$ by two, a contradiction. $\square$

Since there are only a finite number of edges (and simple closed curves) in $\Gamma$, there are only a finite number of isotopy types of simple closed curves in $\eta(\Gamma) \subset \partial B^3$ which can arise as components of $F \cap \partial B^3$. Thus any normal surface can be described completely by saying how many of each of the finite number of possible types occur. Normal surfaces are useful in constructing algorithms because the decisions made in creating a normal surface are essentially finite. The theory is even more powerful than these considerations suggest, since the operation of adding the numbers that classify two surfaces has geometric content. For a full appreciation, it is helpful to consider the very specific type of Heegaard splitting that comes from a triangulation.

8.2. **Special case: Normal surfaces in a triangulation.** Let $M$ be a closed triangulated 3–manifold with a fixed triangulation $\mathcal{T}$. Let $T^i$ denote the $i$–skeleton of $\mathcal{T}$. We will consider what it means for a surface to be normal in the induced Heegaard splitting $H_1 \cup_S H_2$ where, in contrast to 2.1, $H_2$ is a neighborhood of $T^1$ and $H_1$ is a neighborhood of the dual 1-skeleton.

Suppose $F$ is a closed surface in $M$. Then the requirement in 8.1 that $F$ be disjoint from the 0-handles in $H_2$ and intersect the cores of the 1-handles of $H_2$ and $\partial B^3$ transversally here translates to the requirement that $F$ be in general position with respect to the triangulation $\mathcal{T}$. The weight of $F$ is just the number $|F \cap T^1|$. Since $H_1$ is a neighborhood of the dual complex to the triangulation, the 2-simplices of $\mathcal{T}$ are a collection $\Delta_1$ of meridians for $H_1$. The 3-simplices of $\mathcal{T}$ are the balls that are produced when $H_1$ is cut up along $\Delta_1$. The graph $\Gamma$ appears on the boundary of each 3-simplex $\tau$ as the dual (tetrahedral) graph to the 1-skeleton $\tau_1$ of the tetrahedron $\tau$.

Consider how a normal surface intersects the boundary of $\tau$. A single component $c$ of $F \cap \tau_2$ can run only once along any edge of $\Gamma$, or, put another way, $c$ can cross an edge of $\tau_1$ at most once. In particular $c$ meets each face of $\tau$ in a single spanning arc (i. e. an arc whose ends lie on different sides of the triangular face). It follows immediately that a tetrahedron has up to normal isotopy precisely seven curve types. (See fig. 11.)There are four curve types with three sides and three curve types with four sides.

If $\alpha$ is a curve type in $\tau$, and p is a point in the interior of $\tau$, then the cone p*$\alpha$ of $\alpha$ to $p$ is called a *disk type* of $\tau$. Hence a tetrahedron has up to normal isotopy precisely seven disk types. We conclude that $F \subset M$ is a normal surface if and only if $F$ intersects each tetrahedron of $\mathcal{T}$ in a (necessarily pairwise disjoint) collection of these disktypes.

Thus a normal surface is determined by the number of each curve type in which it meets the boundaries of the various tetrahedra. That is, if $\mathcal{C}_1, ..., \mathcal{C}_n$ is an ordering of the curve types, then the surface $F$ determines (and is determined by) an $n$–tuple $(x_1, ..., x_n)$, where $x_i$ denotes the number of representatives of $\mathcal{C}_i$ which $F$ induces in the tetrahedra of $\mathcal{T}$.

Conversely, if we start with an $n$–tuple of non-negative integers, then we can construct a normal surface in $M$ corresponding to this $n$–tuple if it satisfies the following constraints:

1. We can't have two 4–sided disks from distinct normal isotopy classes in the same tetrahedron (for they necessarily intersect).

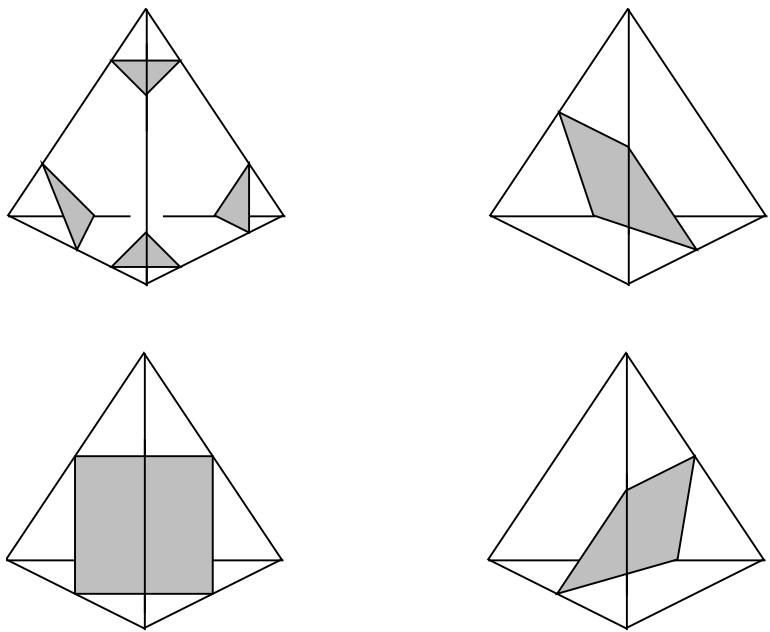

Figure 11.

2. Edges of disktypes on corresponding faces of incident tetrahedra have to match. Namely, if $F$ intersects one face of a tetrahedron in $p$ representatives of a certain arc type, then $F$ also has to intersect the corresponding face of the incident tetrahedron in $p$ representatives of the same arc type.

We now explain "geometric" addition of normal surfaces. A normal surface $F$ in $M$ is *straight* if it satisfies these conditions

1. For any 2−simplex $\sigma$ in $T^2$, $\sigma \cap F$ consists only of straight spanning arcs (called *chords*).
2. In each tetrahedron $\tau$ any 3−sided disk in $\tau \cap F$ is the triangle given by the convex hull of its vertices.
3. Any 4−sided disk in $\tau \cap F$ is the cone to the barycenter of its four vertices.

Clearly any normal surface can be isotoped to be straight. Now consider how two straight normal surfaces $F_1$ and $F_2$ intersect. First move them slightly so that $F_1 \cap F_2 \cap T^1 = \emptyset$ and so that no barycenter of a 4-sided disk in $F_2$ lies in $F_1$ (and vice versa). Then

**Lemma 8.4.** *In each tetrahedron $\tau$, $F_1 \cap F_2$ consists of proper arcs, each of which has its ends on distinct 2−simplices. Each end is a point* *in a 2-simplex $\sigma \prec \tau$ where a chord of $F_1 \cap \sigma$ and a chord of $F_2 \cap \sigma$ intersect.* □

Consider how chords in a 2−simplex $\sigma$ can intersect. Let $p$ be the intersection point. There is a unique way to remove an X neighborhood of $p$ and rejoin the endpoints of the X by two disjoint arcs so that the result gives two spanning arcs in $\sigma$. This process is called a *regular exchange* at $p$.

Now consider extending this regular exchange along an arc component $C$ of $F_1 \cap F_2$ inside a tetrahedron. That is, given two straight disks in a tetrahedron which intersect along an arc $C$, try to remove a neighborhood of $C$ from both $F_1$ and $F_2$ and reattach the sides so that the result is a regular exchange at the ends of $C$. It is easy to see that this is possible, *unless* the disk types are distinct and both 4−sided.

We say that normal surfaces $F_1$ and $F_2$ are *compatible* if, in each tetrahedron, the four-sided curve types of $F_1$ and $F_2$ (if any) are the same. If $F_1$ and $F_2$ are compatible then, after they are straightened, we have seen that in a neighborhood of each curve of $F_1 \cap F_2$ it is possible to perform a regular exchange to eliminate the curve of intersection. The result of this operation on all intersection curves is a normal surface called the *geometric sum* of $F_1$ and $F_2$. Denote this surface by $F_1 + F_2$.

There are several interesting properties which are additive with respect to the geometric sum operation.

If $F_1$ and $F_2$ are compatible normal surfaces, then $F_1 + F_2$ is defined and

1. $\chi(F_1 + F_2) = \chi(F_1) + \chi(F_2)$, where $\chi$ is Euler characteristic.
2. If $F_1$ corresponds to $(x_1, ..., x_n)$ and $F_2$ corresponds to $(y_1, ..., y_n)$, then $F_1 + F_2$ corresponds to $(x_1 + y_1, ..., x_n + y_n)$
3. $w(F_1+F_2) = w(F_1)+w(F_2)$, where $w(F) =$ weight of $F = |F \cap T^1|$.

Now it is easy to see that the solution set of a system of integral equations in the positive orthant is generated under addition by a finite number of "fundamental" solutions which can be found algorithmically (see e. g. [Hm, Chapter 8]). Exploiting the properties of the geometric sum listed above, it's often possible to show that if any surface with a compression preserved property appears in $M$ then one with this property appears among the fundamental surfaces. If it can then be checked whether each of the fundamental surfaces has the desired property, the result is an algorithm to decide if $M$ contains a surface with the desired property. So, for example, there is an algorithm to detect the presence of a reducing sphere, and an algorithm to detect the presence of an injective surface (see [JO], [BS]). Part of this problem requires recognizing if a 2-sphere is a reducing sphere or whether it bounds a 3-ball (see [Th]). This is a difficult problem in its own right, and one that requires a new idea - that of an "almost normal surface". Such a surface is normal, except in a single tetrahedron whose boundary it intersects in an octagon. (See fig. 12) This leads us into an area of very active research. For example, see [St] for a discussion of how strongly irreducible splitting surfaces can be put in almost normal position and see [Ru] for a provocative discussion of other algorithms which may be useful and which make use of almost normal surfaces.

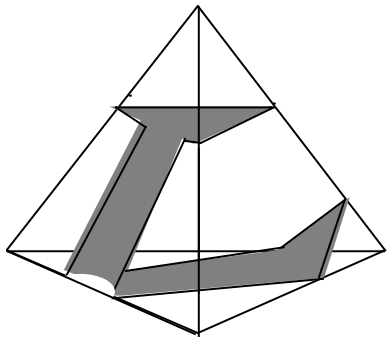

Figure 12.

Martin Scharlemann, Mathematics Department, University of California, Santa Barbara, CA USA

*E-mail address*: mgscharl@math.ucsb.edu